\documentclass[a4paper,10pt]{article}

\usepackage[utf8]{inputenc}
\usepackage{amsmath}
\usepackage{amssymb}
\usepackage [hmargin={3cm,3cm},vmargin={2.4cm,2.6cm}]{geometry}
\usepackage{hyperref}

\newtheorem{thm}{Theorem}[section]

\newtheorem{prop}{Proposition}[section]

\newtheorem{lem}{Lemma}[section]

\newtheorem{rem}{Remark}[section]

\newtheorem{ex}{Exemple}[section]

\newtheorem{cor}{Corollary}[section]

\newtheorem{defin}{Definition}[section]

\newcommand{\R}{\mathbb{R}}

\newcommand{\cqfd}{\begin{flushright} 
                    $\Box$
                   \end{flushright}}

\title{On the finiteness of the Morse index for Schrödinger operators}

\author{Baptiste Devyver\footnote{Laboratoire de Math\'{e}matiques Jean Leray, Universit\'{e} de Nantes; 2 rue de la Houssinière BP 92208, 44322 Nantes cedex 03, France; email: baptiste.devyver@univ-nantes.fr}}

\begin{document}

\maketitle

\begin{keywords}
Spectrum of Schrödinger operators, Morse index, Bochner technique, parabolicity, Doob transform.
\end{keywords}
\textit{MSC Classification: 35P, 58J, 31E}

\begin{abstract}

We show that if $M$ is a complete Riemannian manifold and $H=\Delta+V$ is a Schrödinger operator, then the existence of a positive solution of $Hu=0$ outside a compact set is equivalent to the finiteness of the Morse index of $H$. 

\end{abstract}

\section{Introduction}

A classical and important tool in differential geometry is the so-called Bochner technique. Since its introduction by Bochner, this technique has undergone a huge number of refinements; because it motivates this article, we want to recall a quite general setting to which it applies (see \cite{Pigola-Rigoli-Setti} for details). Consider a Riemannian manifold $M$ and a Riemannian vector bundle $E$ over it, which carries a compatible metric connection $D$. We assume that there is a geometric Laplacian $\vec{\Delta}$ acting on sections of $E$, which can be related to the ``rough Laplacian'' $\bar{\Delta}=-Tr(D^2)$ by the formula:

$$\vec{\Delta}=\bar{\Delta}+\mathcal{R},$$
where $\mathcal{R}$ is a symetric endomorphism in each fiber of $E$. Classical examples of such situations are p-differential forms for the Hodge Laplacian, and spinors for the Dirac Laplacian. Denote by $\mathcal{H}(E)$ the set of section $\xi$ of $E$ such that $\vec{\Delta}\xi=0$. Then, defining $V(x)$ to be the lower eigenvalue of $\mathcal{R}$ and using Kato's inequality, we have:

$$\Delta|\xi|+V|\xi|\leq 0,\,\forall \xi\in\mathcal{H}(E),$$
where $\Delta$ is the Laplacian on $M$ (with the convention that it is a positive operator). Typically, one wants to show that the space of harmonic sections for $\vec{\Delta}$, satisfying some integrability conditions (for exemple, being in $L^2$), has finite dimension. For the $L^2$ integrability conditions, Theorem 5.1 in \cite{Pigola-Rigoli-Setti} asserts in particular that this is the case if we can find a positive function $\varphi$, solution of the equation $\Delta\varphi+V\varphi=0$ outside a compact set. Therefore, in this case the question reduces to give conditions on the potential $V$ such that we can find such a solution. This has a link with the spectrum of the Schrödinger operator $\Delta+V$, as the following Lemma shows (which we extract from \cite{Pigola-Rigoli-Setti}, although it is originally due to Moss and Piepenbrink \cite{Moss-Piepenbrink} and Fischer-Colbrie and Schoen \cite{Fischer-Schoen}):

\begin{lem}\label{Fischer-Schoen}
 
Let $M$ be a Riemannian manifold, $\Omega\subset M$ be a smooth open set, and $V\in L^\infty_{loc}$. Denote by $H_\Omega$ the Schrödinger operator $H:=\Delta+V$ on $\Omega$, with Dirichlet boundary conditions, and assume it is bounded from below. We identify it with its Friedrichs extension, which is self-adjoint. Then the following are equivalent:

\begin{enumerate}
 \item There exists $\varphi\in W^{1,2}_{loc}$ positive solution of
 
 $$H\varphi=0\,\,on\,\,\Omega.$$  
 \item There exists $\varphi\in W^{1,2}_{loc}$ positive such that
 
 $$H\varphi\geq0\,\,on\,\,\Omega.$$
 \item $\lambda_1(H_\Omega)\geq0$, where $\lambda_1$ denotes the infimum of the spectrum.
 
\end{enumerate}

\end{lem}

\begin{defin}
  
We say that a self-adjoint operator has \textbf{finite Morse index} if its essential spectrum $\sigma_{ess}$ is contained in $[0,\infty)$, and if it has only a finite number of negative eigenvalues. In this case, the number (counting multiplicities) of negative eigenvalues is called the \textbf{Morse index} of the operator.

\end{defin}
Notice that for a Schrödinger operator $\Delta+V$, the condition $\sigma_{ess}\subset [0,\infty)$ is satisfied if $V$ tends to $0$ at infinity for example. In \cite{Fischer} (see also \cite{Pigola-Rigoli-Setti}, Chapter 3), Fischer-Colbrie has shown the following Theorem:
 
\begin{thm}\label{Fischer}
 
Let $H:=\Delta+V$ be a Schrödinger operator bounded from below on a complete Riemannian manifold, with $V\in L_{loc}^\infty$, which has finite Morse index. Then there is a positive function $\varphi$, solution of the equation $H\varphi=0$ outside a compact set.

\end{thm}
The proof consists in showing that the finiteness of the Morse index implies that we can find a compact set $K$ such that $\lambda_1(H_{M\setminus K})\geq0$, and applying Lemma (\ref{Fischer-Schoen}).

\begin{rem}
 
We mention that the reason why Fischer-Colbrie and Schoen studied Schrödinger operators was to prove results concerning minimal surfaces. Indeed, for a minimal surface $M$ in a 3-dimensional manifold $N$, we can consider the Schrödinger operator (called the \textbf{stability operator}) $H:=\Delta+S-K+\frac{1}{2}|A|^2$ on $M$, where $S$ is the scalar curvature of $N$, $K$ is the Gaussian curvature of $M$ and $A$ is the second fundamental form of the immersion. $H$ is the linear operator of second variation of the local area functional on the surface $M$; since $M$ is minimal, every point of $M$ is critical point of the local area functional, and saying that $M$ is stable means that up to second order, the deformations of $M$ make the area increase. The spectral properties of $H$, like the fact that $H$ has a finite number of negative eigenvalues, have consequences first for the geometry of $M$, and then for the topology of $N$.

\end{rem}
But it could well happen \textit{a priori} that the existence of a positive solution of $H\varphi=0$ outside a compact set can be true in some cases where there is an infinite number of negative eigenvalues for $H$, converging to zero. In this paper, we show that this cannot happen. Our main result is the following

\begin{thm}\label{main result}

Let $M$ be a complete Riemannian manifold. Let $V\in L^\infty_{loc}$, and denote $H=\Delta+V$ the corresponding Schrödinger operator, which we assume to be bounded from below.
Then the following assertions are equivalent:

\begin{enumerate}
 \item $H$ has finite Morse index.
 
 \item There exists a positive smooth function $\varphi$ in $W^{1,2}_{loc}$ which satisfies $H\varphi=0$ outside a compact set.
 
\end{enumerate}
Furthermore, in this case $Ker_{L^2}(H)$ is finite dimensional.

\end{thm}
We first want to make several remarks concerning this result:

\begin{rem}

The hypothesis that $H$ is bounded from below is to ensure that defined on $C_0^\infty(M)$ --the set of compactly supported smooth functions--, it is an essentially self-adjoint operator.

\end{rem}

\begin{rem}

Theorem (\ref{main result}) also holds for more general operators $H$: the proof will show that in fact it holds for $H$ of the form $(\Delta_{\mu}+W)+V$, where $\mu$ is a $C^1$ positive function and $W\geq0$.  

\end{rem}
There will be two ingredients in the proof. First we consider the operator $L=\varphi^{-1}H\varphi$: this way to transform $H$ is called the \textit{Doob transform} associated to $\varphi$. Then $L$ is also a Schrödinger operator:

$$L=\Delta_{\varphi^2}+\tilde{V},$$
but we will see that its potential $\tilde{V}$ is compactly supported. Thus, Theorem (\ref{main result}) will be a consequence of the following general result, which is of independant interest:

\begin{thm}\label{potentiel support compact}

Let $L$ be an operator of the type: $L=\Delta_{\mu}+W$ with $W\geq0$, where $\Delta_\mu=-\frac{1}{\mu}div (\mu grad)$ with $\mu$ a $C^1$ positive function is a weighted Laplacian. Let $V$ be a compactly supported potential in $L^p$ for a $\frac{n}{2}<p\leq \infty$.\\
Then $\sup\{dim (F) : F\subset C_0^\infty\hbox{ and }q|_F\leq0\}$, where $q$ is the quadratic form associated to $L+V$, is finite. Furthermore, $Ker_{L^2}(L+V):=\{\varphi\in L^2 : (L+V)\varphi=0\}$ has finite dimension.

\end{thm}
Roughly, this Theorem relies on two principles: first, following an idea that goes back to Birman and Schwinger (see for exemple \cite{Reed-Simon4}, p.98-99), we will bound $\sup\{dim F : F\subset C_0^\infty\hbox{ and }q|_F\leq0\}$ by the number of eigenvalues of $L^{-1/2}(-V)L^{-1/2}$ which are greater or equal to $1$. The second idea, which comes from \cite{Carron1}, Proposition 1.2, is that Sobolev inequalities, and more generally non-parabolicity of $(M,g)$ have functional consequences for the operator $\Delta^{-1/2}V\Delta^{-1/2}$: in the case where $V$ has compact support, this operator is compact if $(M,g)$ is non-parabolic. We will extend this to our case of interest, i.e. to $L^{-1/2}(-V)L^{-1/2}$ when $M$ is non-parabolic for $L$ (details of the meaning of this are given in the next two sections). Finally, we will use a trick to deal with the case where $L$ is parabolic.

\begin{rem}
 
As we have already pointed out, the statement concerning the finiteness of the dimension of $Ker_{L^2}(L+V)$ could also be obtained for $L=\Delta$, under the hypothesis that $V$ is continuous, by applying Theorem 5.1 in \cite{Pigola-Rigoli-Setti} (which is much more general than that). Our proof is different, and it has the advantage to be a fairly direct consequence of the results related to the non-parabolicity of $L$ that we prove in section 3.

\end{rem}

\begin{rem}

If $M$ satisfies a Sobolev inequality :

$$||f||_{\frac{2\nu}{\nu-2}}\leq C||\nabla f||_2,\,\forall f\in C_0^\infty(M),$$
then the Cwickel-Lieb-Rosenbljum bound tells us that $N_-(V)$, the cardinal of the negative spectrum of $\Delta+V$, satisfies :

$$N_-(V)\leq C\int_M(V_-)^{\nu/2},$$
where $V_-=-\inf(V,0)$ is the non-positive part of $V$ (see \cite{Saloff-Coste1}).

\end{rem}
\noindent The article is organised as follows: in the first two sections, we develop the material we will need to prove our results. In the first part, we investigate the notion of parabolicity for an operator $L$ of the type : $L=\Delta_\mu+W$. In the second one, we describe functionnal consequences of non-parabolicity for the operator $L^{-1/2}$. In the third part, we prove the two results cited above, in the case of a smooth potential for the proof of Theorem (\ref{main result}). In a fourth one, we weaken the regularity assumptions needed on $V$, and in a fifth one, we present two alternative approaches of the proof of Theorem (\ref{main result}).

\section{On the parabolicity of a manifold}

In this section, we recall the notion of parabolicity. References for this section are \cite{Ancona} and \cite{Grigor'yan1}.\\
\textit{Notations}: Throughout this paper, $(M,g)$ denotes a complete Riemannian manifold, $dx$ is the Riemannian measure on $M$ and $C_0^\infty(M)$ (or $C_0^\infty$ for short) is the set of compactly supported, smooth functions on $M$.\\
We consider on $M$ an operator $L$ of the type $L=\Delta_\mu+W$, $W$ non-negative. It is a well-known fact that $L$ is a positive self-adjoint operator on  $L^2(M,\mu dx)$, associated to the closable quadratic form:

$$q(u)=\int_M(|du|^2 +Wu^2)\mu dx .$$

Particular example of such an operator are:

\begin{enumerate}
 \item The natural Laplacian $\Delta=-div\, grad$, where $div X= \sum_jg^{-1/2}\partial_j(g^{1/2}X^j)$ in a coordinate system.
 
 \item The $\mu$-Laplacian $\Delta_\mu=-\frac{1}{\mu}div(\mu\, grad)$, where $\mu$ is a positive, smooth function on $M$.
 
 \item Schrödinger operators $H=\Delta+W$, with $W\in L^1_{loc}$ a non-negative function.
\end{enumerate}
Return to the general case. We have a Green-type formula: 

\begin{prop}\label{green}

If $u$ and $v$ are elements of $C_0^\infty$,

$$\int_M (u Lv) \mu dx=\int_M (\langle du,dv\rangle +Wuv )\mu dx$$
  
\end{prop}
\textit{Notation:} we denote by $d\nu$ the measure $\mu dx$.

\begin{rem}
 
The restriction $W$ non-negative is to ensure that $L$ satisfies the maximum principle.

\end{rem}
Given $\Omega\subset M$ an open, regular, relatively compact set, let $L_\Omega$ be the self-adjoint operator associated to the restriction of the quadratic form $q$ to the Sobolev space $W_0^{1,2}(\Omega,d\nu)$ (i.e. with Dirichlet conditions). We can consider the Green kernel $G_\Omega$ of $L$ on $\Omega$ with Dirichlet conditions, extended by zero outside $\Omega\times\Omega$; it enjoys the following the properties:

\begin{enumerate}
 \item $G_\Omega\geq0$,

\item  $G_\Omega$ is finite off the diagonal.

\item $G\mid_{\partial(\Omega\times\Omega)}=0$,

\item For all $f\in L^2(\Omega,d\nu)$, $g:=G_\Omega f$ (where $G_\Omega f(x)=\int_\Omega G_\Omega(x,y)f(y) d\nu(y)$) satisfies : 

$$g\in Dom(L_\Omega), \hbox{ and }L_\Omega g=f$$

\end{enumerate}
It is a consequence of the maximum principle that $G_\Omega$ is non-decreasing with respect to $\Omega$:

$$\hbox{if }\Omega_1\subset\Omega_2,\,G_{\Omega_2}\geq G_{\Omega_1},$$
so we can define a pointwise limit:

$$G(x,y):=\lim_{\Omega\rightarrow M}G_\Omega(x,y)\hbox{ for all }x\neq y$$

\begin{defin}

We say that $(M,g)$ is \textbf{non-parabolic} for $L$ if $G(x,y)<\infty$ for a certain $(x,y)$.
 
\end{defin}
As a consequence of the Harnack inequality, we see that this is equivalent to $G(x,y)<\infty$ for all $x\neq y$ (for an account on the notion of parabolicity for the usual Laplacian and a proof of this fact, see the survey of Grigor'yan \cite{Grigor'yan1}). There is a caracterisation of non-parabolicity in term of the ``Dirichlet form'' q of $L$, which we will make constant use of (for a proof and references, see \cite{Ancona}, p.46-47):

\begin{thm}\label{non-parabolicity}
 
The following statements are equivalent: 
\begin{enumerate}
\item $(M,g)$ is non-parabolic for $L$.

\item There exists an open, relatively compact subset $\Omega$ of $M$ and a constant $C\geq0$ such that for all $f\in C_0^\infty(M)$,

$$\int_\Omega f^2d\nu\leq Cq(f)$$

\item Property (2) is true for all open, relatively compact subset $\Omega$ of $M$.
\end{enumerate}

\end{thm}

\begin{cor}\label{non-parabolicity of L+W,W>c>0}
 
If there is an $\epsilon>0$ such that  $W>\epsilon$ on an open set, then $(M,g)$ is non-parabolic for $\Delta_\mu+W$.

\end{cor}

\begin{ex}
 
$\mathbb{R}^n$ is non-parabolic for $\Delta$ if and only if $n>2$. More generally, any complete Riemannian manifold satisfying a Sobolev inequality of index $n>2$:

$$||f||_{\frac{2n}{n-2}}\leq C||\nabla f||_2,\,\forall f\in C_0^\infty$$
is non-parabolic for $\Delta$ (this is an easy consequence of Theorem (\ref{non-parabolicity})). Example of such a manifold other than $\mathbb{R}^n$ is the connected sum of two copies of $\mathbb{R}^n$, for $n>2$.

\end{ex}

\section{Consequences for $L^{-1/2}$}\label{consequence parabolicite}

In this section, we consider as before an operator $L=\Delta_\mu+W$, $W$ non-negative, which is non-parabolic, and we review some functional properties of the operator $L^{-1/2}$ that come from the non-parabolicity of $L$. We keep the notations of section 2. We will define an operator $L^{-1/2}$ by two different means. Finally, we will have to show that these definitions are consistent, in that they agree in a suitable sense.

\begin{defin}
 
We define two unbounded operators $L_s^{1/2}$ and $L_s^{-1/2}$ by the functionnal calculus: if $f$ is a Borel function on $\R$, we can define

$$f(L):=\int_0^\infty f(\lambda)dP_\lambda,$$
where $dP_\lambda$ is the projection-valued measure associated to the self-adjoint operator $L$ (see (\cite{Reed-Simon1})).\\
Then $L_s^{1/2}:=f(L)$ with $f(x)=x^{1/2}$, and $L_s^{1/2}:=g(L)$ with $g(x)=x^{-1/2}$.

\end{defin}

\begin{rem}
 
\begin{enumerate}

\item Since $(M,g)$ is non-parabolic for $L$, the functionnal inequalities of Theorem (\ref{non-parabolicity}) imply that $Ker_{L^2}L=\{0\}$, for if $u\in Ker_{L^2}L$ we have $q(u)=0$ by definition of $q$. Therefore, $P_{\{0\}}=0$ and we can indeed take $g(x)=x^{-1/2}$ in the above definition, even if $g$ is not defined in $0$.

\item The ``s'' index stands for ``spectral''.

\item By construction, $\mathcal{D}(L_s^{1/2})=\mathcal{D}(q)$ (where $\mathcal{D}$ denotes the domain).

\end{enumerate}

\end{rem}
The non-parabolicity of $(M,g)$ for $L$ allows us to consider an alternative definition of $L^{-1/2}$, which we describe now. Let $H_0^1$ be the closure of $C_0^\infty(M)$ for the norm $N(u)=|| L_s^{1/2}u||_2=(\int (|du|^2+Wu^2)d\nu)^{1/2}$. It is a Hilbert space, and we have the following paraphrase of the implication $(1)\Rightarrow (3)$ of Theorem (\ref{non-parabolicity}), which allows us to see $H_0^1$ as a functions space:

\begin{prop}\label{Sobolev space definition}

If $M$ is non-parabolic for $L$, then the injection $C_0^\infty(M)\hookrightarrow W_{loc}^{1,2}(M,d\nu)$ extends continuously to:

$$H_0^1 \hookrightarrow W^{1,2}_{loc}(M),$$
that is : for all $U$ open, relatively compact set, the restriction to $U$ of elements of $H_0^1$ belong to $W^{1,2}(U)$, and there exists a constant $C_U$ such that 

$$\big|\big| f|_U\big|\big|_{W^{1,2}(U)}\leq C_U||f||_{H_0^1},\,\forall f\in H_0^1,$$
or equivalently:

$$\int_U f^2\leq C_U||f||^2_{H_0^1},\,\forall f\in H_0^1.$$

\end{prop}
We then define:

\begin{defin}
 
The operator $L_s^{1/2}$, restricted to $C_0^\infty(M)$, extends to an isometry:

$$L_a^{1/2} : H_0^1\longrightarrow L^2(M, d\nu)$$

\end{defin}
The following Proposition tells us that these two operators $L_s^{1/2}$ and $L_a^{1/2}$ are in fact equal:

\begin{prop}\label{coincidence operateur spectral et absolu}

\begin{enumerate}
\item $C_0^\infty(M)$ is a core for $L_s^{1/2}$.

\item$\mathcal{D}(L_s^{1/2})=H_0^1\cap L^2$, and the restrictions to
$H_0^1\cap L^2$ of the two operators $L_a^{1/2}$ and $L_s^{1/2}$
are equal.
\end{enumerate}

\end{prop}
\textit{Proof:}\\
(1): Let $A$ be the restriction of $L_s^{1/2}$ to $C_0^\infty(M)$. We have to show that

$$Im(A\pm i)^\perp=\{0\}$$
Let $f\in Im(A+i)^\perp$. Then for all $g\in C_0^\infty(M)$, 

$$\langle f,(A+i) g\rangle=0.$$
We can write $f=(L_s^{1/2}+i)u$, where $u\in L^2$, since $L_s^{1/2}$ is self-adjoint. Then if $g\in C_0^\infty(M)$, 

$$0=\langle f,(A+i)g\rangle=\langle (L_s^{1/2}+i)u,(L_s^{1/2}+i)g\rangle=\langle u,(L_s+1)g\rangle,$$ 
by the Spectral Theorem. But $C_0^\infty(M)$ is a core for $L_s+1$, therefore $u\in \mathcal{D}(L_s+1)$ and $(L_s+1)u=0$. Since $L_s\geq 0$, $-1$ does not belong to the spectrum of $L_s$, and we conclude that $u=0$, then $f=0$.\\
The proof for $A-i$ is similar.

\vskip3mm

\noindent (2): Define a quadratic form $Q$ on $C_0^\infty$ by 

$$Q(f)=||f||^2_2+\langle Lf,f\rangle=||f||^2_2+||L_s^{1/2}f||_2^2,$$
and a quadratic form $\bar{Q}$ on $H_0^1\cap L^2$ by 

$$\bar{Q}(f)=||f||_2^2+||L_a^{1/2}f||_2^2= ||f||_2^2+||f||_{H_0^1}^2.$$
By the consequence of non-parabolicity given in Proposition (\ref{Sobolev space definition}), $\bar{Q}$ is closed. It is thus a closed extension of $Q$, which yields a self-adjoint operator $S$ such that $\mathcal{D}(S^{1/2})=\mathcal{D}(\bar{Q})$ and for all $f\in \mathcal{D}(\bar{Q}),$

$$\bar{Q}(f)=\langle Sf,f\rangle=||S^{1/2}f||_2^2.$$
But since $L_s$ is essentially self-adjoint, it has a unique self-adjoint extension, and so we get that $S=L+1$. Using the fact that

$$\sqrt{L_s}\leq \sqrt{L_s+1}\leq \sqrt{L_s}+1,$$
we get that $\mathcal{D}(\sqrt{L_s+1})=\mathcal{D}(\sqrt{L_s})$, and then $\mathcal{D}(\bar{Q})=\mathcal{D}(L_s^{1/2})$. By the first part of Proposition (\ref{coincidence operateur spectral et absolu}), $\sqrt{L_s+1}|_{C_0^\infty(M)}$ is essentially self-adjoint, and since

$$\bar{Q}(f)=\big|\big|\sqrt{L_s+1}f\big|\big|_2^2,\,\forall f\in \mathcal{D}(\bar{Q}),$$ 
we conclude that $C_0^\infty(M)$ is dense in $\mathcal{D}(\bar{Q})=H_0^1\cap L^2$ for the norm given by $\sqrt{\bar{Q}}$.\\
Since $L_a^{1/2}$ and $L_s^{1/2}$ coincide on $C_0^\infty(M)$, by a limit argument they also coincide on $H_0^1\cap L^2$.

\cqfd
From Proposition (\ref{coincidence operateur spectral et absolu}), we can deduce the following Lemma:

\begin{lem}\label{absolute}
 
If $u$ and $v$ belong to $H_0^1\cap L^2$, then

$$\langle L_a^{1/2}u,v\rangle=\langle u,L_a^{1/2}v\rangle$$

\end{lem}
\textit{Proof:}\\
It is a consequence of the facts that $L_a^{1/2}$ and $L_s^{1/2}$ coincide on $H_0^1\cap L^2$ (by Proposition (\ref{coincidence operateur spectral et absolu})), and that $L_s^{1/2}$ is self-adjoint.

\cqfd

\begin{prop}\label{H_0^1 L^2}

$L_a^{1/2} : H_0^1\rightarrow L^2$ is an isomorphism.

\end{prop}
\textit{Proof of Proposition (\ref{H_0^1 L^2}):}\\
$L_a^{1/2}$ is the unique continuous extension of the isometry $L_s^{1/2}:
C_0^\infty(M)\rightarrow L^2$, so it is also an isometry, hence
injective.\\
To prove that it is onto, since the image of $L_a^{1/2}$ is closed
by the fact that it is an isometry, it is enough to prove that
$(Im\, L_a^{1/2})^{\perp}=\{0\}$. So let $ w\in (Im\, L_a^{1/2})^{\perp}\subset L^2$. Then for all $u\in C_0^\infty(M)$,

$$\langle w, L_s^{1/2}u\rangle =0.$$
Since $C_0^\infty(M)$ is a core for $L_s^{1/2}$, we obtain:

$$w\in Dom(L_s^{1/2})\hbox{ and }L_s^{1/2}w=0.$$
%%But since $L_s$ is essentially self-adjoint over $C_0^\infty$ by
%%unicity of the self-adjoint extension, $L_s$ is the Friedrichs
%%extension of the form $q$ restricted to $C_0^\infty$. It implies
%%that $Dom(q)$ is the closure of $C_0^\infty$ for the norm
%%$||u||_2+||L_s^{1/2}u||_2$, and therefore is included in
%%$H_0^1\cap L^2(M,d\nu)$. By definiton, $Dom(q)=Dom(L_s^{1/2})$, so
%%$Dom(L_s^{1/2})\subset H_0^1$.
We deduce by Lemma (\ref{absolute}) that $w\in H_0^1$ and $L_a^{1/2}w=0$. Since
$L_a^{1/2}$ is injective, $w=0$ in $H_0^1$ and then in $L^2$ by
Theorem (\ref{non-parabolicity}).

\cqfd
To sum up, we have defined an operator:

$$L^{-1/2}=L_s^{-1/2}=L_a^{-1/2},$$
which enjoys the following properties:

\begin{enumerate}

\item It is a bijective isometry from $L^2$ to $H_0^1$.

\item As a non-bounded operator, it has domain $L^{1/2}(H_0^1\cap L^2)$.

\end{enumerate}
Later, we will look to the operator $L^{-1/2}VL^{-1/2}$, when $V$ is compactly supported (here, we have identified $V$ with the operator ``multiplication by $V$''). To show that it is compact, we will need the following:

\begin{prop}\label{fonctions a support compact dans domaine}
 
Denote by $L_0^2$ the space of compactly supported $L^2$ functions.\\
Then $L_0^2$ is contained in $Dom(L_s^{-1/2})$, and given an open, relatively compact subset $K$, there exists a constant $C_K$ such that for all $v\in L_0^2$ with support included in $K$,

$$||L_s^{-1/2}v||_2\leq C_K||v||_2.$$
More generally, if $n$ denotes the dimension of $M$, for all $n\leq p\leq \infty$, there is a constant $C(p,K)$ such that for all $v\in L_0^{\frac{2p}{p+2}}$ with support included in $K$, 

$$||L_s^{-1/2}v||_2\leq C(K,p)||v||_{\frac{2p}{p+2}}.$$

\end{prop}
\textit{Proof of Proposition (\ref{fonctions a support compact dans domaine}):}\\
The proof is by duality: we will make use of the following

\begin{lem}\label{duality for functions with compact support}
 
To each $v\in L_0^{\frac{2p}{p+2}}$, we associate the linear form

$$\varphi_v : w\in H_0^1\cap L^{\frac{2p}{p-2}}\hookrightarrow \langle v,w\rangle$$
Then for all $n\leq p\leq \infty$, $\varphi_v$ extends uniquely to an element of $(H_0^1)^\prime$, with $||\varphi_v||\leq C(p,K)||v||_{\frac{2p}{p+2}}$, where $K$ is an open, relatively compact subset containing the support of $v$.

\end{lem}
The proof of Lemma (\ref{duality for functions with compact support}) is a consequence of non-parabolicity of $L$, and we prove it after.\\
We then conclude the proof of Proposition (\ref{fonctions a support compact dans domaine}). Let $u\in H_0^1$ such that $L_a^{1/2}u\in L_0^2$. We have to show that $u\in L^2$. Let $v\in L_0^2$ defined by $v:=L_a^{1/2}u$. By Lemma (\ref{duality for functions with compact support}), there exists $h\in H_0^1$ such that $\varphi_v=\langle h,.\rangle_{H_0^1}=\langle L^{1/2}h,L^{1/2} .\rangle$. We define $f=L^{1/2}h$, and since $\mathcal{D}(L_s^{1/2})=H_0^1\cap L^2$, we obtain that $f\in \mathcal{D}((L_s^{1/2})^*)$ with $(L_s^{1/2})^*f=v$. But $L_s^{1/2}$ is self-adjoint, so $f\in H_0^1\cap L^2$ and $L_a^{1/2}f=L_s^{1/2}f=v=L_a^{1/2}u$. $L_a^{1/2}$ being injective, $u=f\in L^2$.\\
For the inequality on the norm, we remark that $||u||_{L^2}=||f||_{L^2}=||h||_{H_0^1}=||\varphi_{L^{1/2}u}||$, and use Lemma (\ref{duality for functions with compact support}).

\cqfd
\textit{Proof of Lemma (\ref{duality for functions with compact support}):}\\
Let $K$ be an open, relatively compact set containing the support of $v$ in its interior. Let $\rho\in C_0^\infty(M)$ such that $\rho=1$ on the support of $v$, and $\rho=0$ outside $K$. If $w\in C_0^\infty(M)$,

$$|\varphi_v(w)|\leq ||v||_{\frac{2p}{p+2}}||\rho w||_{L^{\frac{2p}{p-2}}(K)}.$$
We first treat the case $p=\infty$, i.e. $\frac{2p}{p-2}=\frac{2p}{p+2}=2$. In this case, we estimate $||\rho w||_2$ by $||\rho||_\infty ||w||_2$. By non-parabolicity, there exists $C_K$ such that 

$$||w||_{L^2(K)}\leq C_K ||w||_{H_0^1},$$
independantly of $w$.\\
Therefore

$$|\varphi_v(w)|\leq C_K||v||_{L^2}||w||_{H_0^1},$$
which proves the result in this case.

\vskip3mm

\noindent For $n\leq p<\infty$, we use the fact that $K$ satisfies the Sobolev inequality:

$$||u||_{\frac{2p}{p-2}}\leq C||\nabla u||_2,\,\forall u\in C_0^\infty(M)\hbox{ such that }supp\, u\subset K.$$
So 

$$||\rho w||_{L^{\frac{2p}{p-2}}(K)}\leq C||\nabla(\rho w)||_2\leq C'\left(||w||_{L^2(K)}+||\nabla w||_{L^2(K)}\right).$$
Since $L=\Delta_\mu+W$ with $W\geq0$ and $\mu$ is bounded from below by a positive constant over $K$, we have, for all $u\in C_0^\infty(M)$: 

$$||\nabla u||_{L^2(K)}\leq C||u||_{H_0^1}.$$ 
We then conclude as before.

\cqfd
We obtain immediately the following:

\begin{cor}\label{boundedness L^-1/2V}
 
For $V\in L_0^p$ with $n\leq p\leq\infty$, the operator $T:= L_s^{-1/2}V$ is bounded on $L^2$.

\end{cor}
\textit{Proof of Corollary (\ref{boundedness L^-1/2V}):}\\
If $K$ is a compact set containing the support of $V$, the operator ``multiplication by $V$'' is bounded from $L^2$ to $L^{\frac{2p}{p+2}}(K)$. We then apply Proposition (\ref{fonctions a support compact dans domaine}).

\cqfd
Furthermore, the non-parabolicity of $M$ yields:

\begin{prop}\label{compactness VL^-1/2}
 
Let $M$ be non-parabolic with respect to $L$. If $V\in L_{loc}^p$ for a $n< p\leq +\infty$, with compact support the operator:

$$VL_a^{-1/2} : L^2\longrightarrow L^2$$
is compact.

\end{prop}
\textit{Proof of Proposition (\ref{compactness VL^-1/2}):}\\
Let $K$ be an open, relatively compact subset of $M$ containing the support of $V$. We can assume that $K$ is smooth. Let $\rho\in C_0^\infty$ such that $\rho|_K=1$. The non-parabolicity criterion of Theorem (\ref{non-parabolicity}) means that:

$$L_a^{-1/2} : L^2\longrightarrow W_{loc}^{1,2}$$
We consider the following compositions:
$$W_{loc}^{1,2}\rightarrow W^{1,2}(K)\hookrightarrow L^{\frac{2p}{p-2}}(K)\rightarrow L^2(K),$$
where the arrow on the left is the multiplication by $\rho$, the one in the middle is the compact Sobolev inclusion (here we use $p>n$), and the one on the right is the multiplication by $V$. The resulting composition is thus compact, and it is in fact equal to the operator ``multiplication by $V$'', sending $W_{loc}^{1,2}$ into $L^2(K)$. Thus we get the result.
\cqfd
%Le fait sur la compacite des injections de Sobolev vient de :
%\begin{prop}
%Soit $T_\infty : W^{1,2}(K)\rightarrow L^2(K)$ l'injection compacte (par Rellich), et $T_n : W^{1,2}(K)\rightarrow L^{\frac{2n}{n-2}}$ l'injection %continue. Alors $T_\infty$ et $T_n$ induisent par interpolation des injections compactes : 
%$$T_p : W^{1,2}(K)\rightarrow L^{\frac{2p}{p-2}}(K),\,\hbox{pour }p>n.$$ 
%\end{prop}
%Preuve :
%Puisque pour $p>n$, $2\leq \frac{2p}{p-2}\leq \frac{2n}{n-2}$, $T_\infty$ et $T_n$ induisent bien une injection $T_p$ par interpolation.
%Si $u_n\in W^{1,2}(K)$ converge faiblement, alors par l'inegalite d'interpolation (consequence de Holder),
%$$||u_n-u_m||_{\frac{2p}{p-2}}\leq ||u_n-u_m||_{frac{2n}{n-2}}^\theta ||u_n-u_m||_2^{1-\theta},$$
%ou $\theta$ verifie :
%$$1-\frac{2}{p}=\theta(1-\frac{2}{n})+\frac{1-\theta}{2}.$$
%Puisque $T_\infty$ est compacte, $u_n$ converge fortement en norme $frac{2n}{n-2}$, et donc par l'inegalite precedente, $(u_n)$ converge en norme $\frac{2p}{p-2}$.$\Box$
Finally, our main result for this section is:

\begin{thm}\label{compact}

Let $V\in L_0^q$ for a $\frac{n}{2}< q\leq +\infty$ be a non-negative, compactly supported potential.\\
Then the operator:

$$L_s^{-1/2}VL_a^{-1/2} : L^2\longrightarrow L^2$$
is self-adjoint, compact.

\end{thm}
\textit{Proof of Theorem (\ref{compact}):}\\
We write:

$$L_s^{-1/2}VL_a^{-1/2}=(L_s^{-1/2}W_1)(W_2L_a^{-1/2}),$$
with $W_1=W_2=V^{1/2}\in L^{p}$ and $p=2q> n$. Let $T_1= L_s^{-1/2}W_1$, and $T_2=W_2L_a^{-1/2}$. By Corollary (\ref{boundedness L^-1/2V}), $T_1 : L^2\rightarrow L^2$ is bounded, and by Proposition (\ref{compactness VL^-1/2})), $T_2: L^2\rightarrow L^2$ is compact. Therefore, $T: =L_s^{-1/2}VL_a^{-1/2}=T_1 T_2$ is compact.\\
To show that $T$ is self-adjoint, we consider first the case where $V\in L_0^\infty$; as before, we decompose $V=W_1W_2$, with $W_1=W_2=V^{1/2}$. It is enough to prove that in this case, $T_1^*=T_2$, i.e. that for all $u,v\in L^2$,

$$\langle L_s^{-1/2}W u,v\rangle =\langle u, W L_a^{-1/2} v\rangle.$$
It is a consequence of a small variation of Lemma (\ref{absolute}): define $f=L_s^{-1/2}W u$ and $g=L_a^{-1/2}v$, then we want to prove:

$$\langle f,L^{1/2}g\rangle=\langle L^{1/2}f,g\rangle$$
Lemma (\ref{absolute}) asserts that it is true if $f,g\in L^2\cap H_0^1$, which happens if $v\in \mathcal{D}:=L^{1/2}(H_0^1\cap L^2)$. Now, $\mathcal{D}$ contains $L^{1/2}(C_0^\infty)$, hence is dense in $L^2$, and we conclude by continuity of $T_1$ and $T_2$.\\
Let us return to the general case. We take an approximation sequence $(V_k)$: $V_k:=\inf(k,V)$. For all $k$, $V_k\in L_0^\infty$ and $V_k\rightarrow V$ in $L^q$-norm; furthermore, the support of $V_k$ is contained in the support of $V$. Define $T_{1,k}:= V_k^{1/2}L_a^{-1/2}$ and $T_{2,k}:=L_s^{-1/2}V_k^{1/2}$. We have $V_k^{1/2}\rightarrow V^{1/2}$ in $L^p$-norm, so by the proof of Proposition (\ref{fonctions a support compact dans domaine}) (resp. by the proof of Proposition (\ref{compactness VL^-1/2})), $T_{1,k}$ (resp. $T_{2,k}$) converges to $T_1$ (resp. to $T_2$) for the strong topology of operators (i.e. $\forall u\in L^2,\,T_{i,k}u\rightarrow T_iu$ in $L^2$). We conclude that the sequence of operators $(L_s^{-1/2}V_kL_a^{-1/2})_k$ converges to $L^{-1/2}_sVL_a^{-1/2}$ for the strong topology of operators. Since each of the $L_s^{-1/2}V_kL_a^{-1/2}$ is self-adjoint, $L^{-1/2}_sVL_a^{-1/2}$ is also self-adjoint.

\cqfd

\section{Main result}

\subsection{A preliminary result}

Now we prove Theorem (\ref{potentiel support compact}). For a potential $V$, define $N_-(V)$ to be the cardinal of $Spec(L+V)\bigcap (-\infty,0)$. We recall two other equivalent definitions of $N_-(V)$, the second one using the fact that $L+V$ is essentially self-adjoint on $C_0^\infty(M)$ (for a proof, see for example \cite{Reed-Simon4}). Let us recall that we denote by $q$ the quadratic form associated to $L+V$.

\begin{prop}\label{nombre valeurs propres negatives}
 
$$\begin{array}{llll}
   
N_-(V)&=&\sup\{dim(F) : F \subset Dom(q)\hbox{ and }q|_F\hbox{ negative definite}\}\\

&=&\sup\{dim(F) : F \subset C_0^\infty(M)\hbox{ and }q|_F\hbox{ negative definite}\}

\end{array}$$

\end{prop}
Remark that the definition 

$$N_-(V)=\sup\{dim(F) : F \subset C_0^\infty\hbox{ and }q|_F\hbox{ negative definite}\}$$
makes sense even if $L+V$ is not essentially self-adjoint. In fact, we will prove that with this definition of $N_-(V)$ and without assuming $L+V$ to be essentially self-adjoint, the existence of a positive solution $\varphi$ of $(L+V)\varphi=0$ outside a compact set implies that $N_-(V)$ is finite.\\
\textit{Proof of theorem (\ref{potentiel support compact}):}\\ 
Since $N_-(V)\leq N_-(-V_-)$, we can assume that $V$ is non-positive. We divide the proof into two steps:\\
\underline{Step 1: case where $L$ is non-parabolic:}\\
For this, we need the following two lemmas:

\begin{lem}\label{valeurs propres}

As in section 2, we denote by $H_0^1$ the space naturally associated to $L$. Let $u\in C_0^\infty(M)$, such that $\langle (L+V)u,u\rangle\leq 0$. Define $v:=L^{1/2}u$.\\
Then 

$$||v||_2^2\leq\langle L^{-1/2}(-V)L^{-1/2}v,v\rangle.$$

\end{lem}
\textit{Proof of Lemma (\ref{valeurs propres}):}\\
The hypothesis is that:

$$\langle Lu,u\rangle\leq\langle (-V)u,u\rangle$$
Using that $\langle Lu,u\rangle=\langle L^{1/2}u,L^{1/2}u\rangle=||v||_2^2$, we get:

$$||v||_2^2\leq\langle (-V)L^{-1/2}v,L^{-1/2}v\rangle,$$
and it remains to prove that:

\begin{equation}\label{e:symmetrie}
 \langle (-V)L^{-1/2}v,L^{-1/2}v\rangle = \langle L^{-1/2}(-V)L^{-1/2}v,v\rangle
\end{equation}
Let $w:=L^{-1/2}(-V)L^{-1/2}v=L^{-1/2}(-V)u$. Equality (\ref{e:symmetrie}) is equivalent to:

$$\langle L^{1/2}w,u\rangle=\langle w,L^{1/2}u\rangle.$$
Now, $u$ belongs to $H_0^1\cap L^2$, and since the operator $L^{-1/2}VL^{-1/2} : L^2\rightarrow L^2$ is bounded, we get $w\in L^2$. In addition, by H\"{o}lder's inequality, $Vu\in L^2$ (since $u\in C_0^\infty$), so $w\in H_0^1$. We conclude by applying Lemma (\ref{absolute}).

\cqfd

\begin{lem}\label{dimension finie}
 
If $S$ is a subspace of $L^2$ such that $S\subset\{v\in L^2 : ||v||_2^2\leq\langle Tv,v\rangle\}$, where $T:=L^{-1/2}(-V)L^{-1/2}$, then the dimension of $S$ is less than the number (counting multiplicities) of eigenvalues of $T$ that are greater than $1$.

\end{lem}
\textit{Proof of lemma (\ref{dimension finie}):} it is an easy consequence of the min-max principle.

\cqfd
\textit{End of the proof of Step 1:}\\
Let $F\subset C_0^\infty$ such that $L+V\leq 0$ on $F$. We have, by definition, $F\subset H_0^1\cap L^2$. Define $S:=L^{1/2}F\subset L^2$. By Lemma (\ref{valeurs propres}), $S\subset\{v\in L^2 : ||v||_2^2\leq\langle Tv,v\rangle\}$, so by Lemma (\ref{dimension finie}), we get that the dimension of $S$ is less than the number of eigenvalues greater than $1$ of $T$. Since $L^{1/2}$ is injective, $dim(F)=dim(S)$, so $dim(F)$ is less than the number of eigenvalues greater than $1$ of $T$. By Theorem (\ref{compact}), $T$ is a self-adjoint compact operator, so the number of its eigenvalues greater than $1$ is finite. Since by Proposition (\ref{nombre valeurs propres negatives}) $N_-(V)$ is less than the number of eigenvalues of $T$ greater than one, $N_-(V)$ is finite, which concludes the first step.\\
\underline{Step 2: general case ($L$ is no more assumed to be non-parabolic):}\\
We write:

$$L+V=(L+\rho)+(V-\rho),$$
where $\rho\in C_0^\infty$ is a non-negative function such that $\rho|_U\geq 1$ for an open set $U$. Define $\tilde{L}:=L+\rho$, and $\tilde{V}:=V-\rho$, so that $\tilde{L}+\tilde{V}=L+V$. By Corollary (\ref{non-parabolicity of L+W,W>c>0}), $\tilde{L}$ is non-parabolic, so we can apply Step 1 to $\tilde{L}+\tilde{V}$, to conclude that it has a finite number of negative eigenvalues.\\\\
It remains to prove the second part of the Theorem, i.e. to prove that $Ker_{L^2}(L+V)$ has finite dimension. As above, it is enough to treat the case where $L$ is non-parabolic. The result is then a consequence of the following:

\begin{lem}\label{noyau}

If $V\in L_0^p$ for some $p>\frac{n}{2}$, then 
 
$Ker_{L^2}(L+V)\subset H_0^1$, and $L^{1/2}Ker_{L^2}(L+V)\subset Ker_{L^2}(I+L^{-1/2}VL^{-1/2})$.

\end{lem}
Given that $L^{1/2} : H^1_0\rightarrow L^2$ is injective, and that $L^{-1/2}VL^{-1/2} : L^2\rightarrow L^2$ is compact by Theorem (\ref{compact}), we obtain that $Ker_{L^2}(L+V)$ is of finite dimension.

\cqfd 
\textit{Proof of Lemma (\ref{noyau}):}\\
The proof is inspired by the Proof of Proposition 1.4 in \cite{Carron1}. Let $\varphi\in Ker_{L^2}(L+V)$. We have $L\varphi=-V\varphi$. We first prove that $\varphi\in H_0^1$. Let $\rho\in C_0^\infty(M)$, then for all $u\in C^\infty(M)$

$$\begin{array}{rcl}
||L^{1/2}(\rho u)||_2^2&=&\langle L(\rho u),\rho u\rangle\\\\
	&=& \int_M \left(\rho^2uLu-\langle \frac{d(\rho^2)}{2},d(u^2)\rangle+\rho u^2 \Delta_\mu\rho\right)d\nu\\\\
	&=& \int_M \left(\rho^2uLu+|d\rho|^2u^2\right)d\nu
\end{array}$$
where we have used integration by parts for the last step. We take $\varphi_k\in C_0^\infty(M)$ such that $\varphi_k\rightarrow \varphi$ in $L^2$ and $L\varphi_k\rightarrow L\varphi$ in $L^2$ (this is possible since $L$ is essentially self-adjoint on $C_0^\infty(M)$). Then, applying the preceding formula, we find that:

$$\lim_{k\rightarrow \infty}||L^{1/2}(\rho \varphi_k)||_2^2= \int_M \left((\rho \varphi)^2 (-V)+|d\rho|^2u^2\right)d\nu,$$
and since $L^2-\lim_{k\rightarrow\infty}\rho\varphi_k=\rho\varphi$ and the quadratic form associated to $L^{1/2}$ is closed, we can let $k\rightarrow\infty$ in the preceeding formula:

$$||L^{1/2}(\rho \varphi)||_2^2= \int_M \left((\rho \varphi)^2 (-V)+|d\rho|^2u^2\right)d\nu.$$
Now we fix a point $o\in M$, and we take a sequence $\rho_k\in C_0^\infty$, such that $\rho_k\equiv 1$ on $B(o,k)$, $\rho_k\equiv 0$ outside $B(o,k+1)$ and $||d\rho_k||_\infty\leq 2$. By applying the formula to $(\rho_k-\rho_l)$, we have for $l\geq k$ such that supp$(V)\subset B(o,k)$,

$$\lim_{k,l\rightarrow\infty}||L^{1/2}((\rho_k-\rho_l) \varphi)||_2^2\leq 4 \lim_{k,l\rightarrow\infty}\int_{M\setminus B(o,k)} u^2d\nu=0,$$
which shows that $L^{1/2}\varphi\in L^2$.

%Let $w\in C_0^\infty$, then by definition,

%$$\langle \varphi,Lw\rangle=-\langle V\varphi,w\rangle.$$
%Define $u:=L^{1/2}w$, which satisfies $u\in L^2\cap H_0^1$ and $-\langle V\varphi,w\rangle=\langle \varphi,Lw\rangle=\langle \varphi, L^{1/2}u\rangle$. We notice that $L^{1/2}C_0^\infty$ is dense in $H_0^1$: this comes from the fact that $L^{1/2}$ is injective, and self-adjoint. Therefore, in order to prove that $\varphi\in H_0^1$, we only have to show that we can find a constant $C$ independant of $w$ such that $|\langle \varphi,L^{1/2}u\rangle|\leq C||u||_2$ (recall that $u=L^{1/2}w$). But the following formula holds:

%$$\langle V\varphi,L^{-1/2}u\rangle =\langle L^{-1/2}V\varphi,u\rangle.$$
%Indeed, by Lemma (\ref{absolute}), it is true for all $\varphi\in C_0^\infty$, $V\in C_0^\infty$, and then recalling by Lemma (\ref{boundedness L^-1/2V}) that $L^{-1/2}V : L^2\rightarrow L^2$ is bounded with norm less that a constant times $||V||_p$, we can take an approximation sequence for $V$ and $\varphi$ and pass to the limit in the equality. Finally, we get:

%$$\langle \varphi,L^{1/2}u\rangle=-\langle L^{-1/2}V\varphi,u\rangle,$$
%and since $L^{-1/2}V : L^2\rightarrow L^2$ is bounded, $\varphi\in H_0^1$.\\
\noindent Given this, we have:

$$\langle \varphi,(L+V)w\rangle=0,\,\forall w\in C_0^\infty.$$
We write $(L+V)w=L^{1/2}(L^{1/2}+L^{-1/2}V)w$; $(L^{1/2}+L^{-1/2}V)w\in L^2$ since $L^{-1/2}V : L^2\rightarrow L^2$ is bounded. Moreover, $L^{1/2}(L^{1/2}+L^{-1/2}V)w=(L+V)w\in L^2$, hence $(L^{1/2}+L^{-1/2}V)w\in \mathcal{D}(L^{1/2})=H_0^1\cap L^2$. Since $\varphi\in H_0^1\cap L^2$, we can apply Lemma (\ref{absolute}) to get:

$$\langle \varphi,(L+V)w\rangle=\langle L^{1/2}\varphi,(L^{1/2}+L^{-1/2}V)w\rangle.$$
Now, $(L^{1/2}+L^{-1/2}V)w=(I+L^{-1/2}VL^{1/2})(L^{1/2}w)$. Define $u=L^{1/2}w\in L^2$, then we have:

$$\langle L^{1/2}\varphi,(I+L^{-1/2}VL^{1/2})u\rangle=0.$$
Since $L^{1/2}C_0^\infty$ is dense in $L^2$, the preceeding equality holds for all $u\in L^2$, and since $(I+L^{-1/2}VL^{1/2})$ is self-adjoint we deduce that $L^{1/2}\varphi\in Ker_{L^2}(I+L^{-1/2}VL^{1/2})$. 

\cqfd

\subsection{Proof of the main result}

The aim of this section is to prove the announced result (Theorem (\ref{main result})), under the supplementary assumption that the potential $V$ is smooth.\\

\noindent \textit{Proof of Theorem (\ref{main result}):}\\
As we have already said, the fact that $N_-(V)<\infty$ implies the existence of a positive solution $\varphi$ of $H\varphi=0$ outside a compact set was proved by Fischer-Colbrie in \cite{Fischer}. The fact that this solution is smooth if $V$ is smooth comes from elliptic regularity. The function $\varphi$ can be smoothly extended to a positive function over $M$. Now we assume the existence of such a solution $\varphi$, and we want to prove that $Card(Spec(H)\cap(-\infty,0])$ is finite. If $u\in L^2$ is an eigenfunction of $H$, i.e. $Hu=\lambda u$ for some $\lambda$, we can write (since $\varphi>0$): 

$$u=v\varphi,$$
then

$$(\varphi^{-1}H\varphi)v=\lambda v$$
Furthermore, if we denote by $d\nu$ the measure $\varphi^2 dx$, we have $v\in L^2(d\nu)$. So we are led to consider the \textit{Doob transform}, which is the following unitary transformation 

\begin{displaymath}
\begin{array}{rll}
 
L^2(d\nu)&\rightarrow& L^2(dx)\\

w&\mapsto& \varphi w

\end{array}
\end{displaymath}
Under this transformation, the operator on $L^2(d\nu)$ associated to $H$ is $L:=\varphi^{-1}H\varphi$. Since the operators $H$ and $L$ are conjugated by a unitary transformation, they have the same spectrum. It turns out that $L$ can be described in another way, thanks to the equation $H\varphi=0$ satisfied by $\varphi$ outside a compact:

\begin{lem}\label{Doob}
 
$$L=\Delta_{\varphi^2}+q,$$
as operators on the distributions, where $q:=\varphi^{-1}H\varphi$ is a compactly supported potential. 

\end{lem}
\textit{Proof of Lemma (\ref{Doob}):}\\
If $v\in C_0^{\infty}(M)$,

$$\begin{array}{llll}
H(\varphi v) &=&\Delta(\varphi v)+V\varphi v \\\\
             &=&(\Delta \varphi)v+\varphi(\Delta v)-2\langle d\varphi,dv\rangle+V\varphi v \\\\
             &=&(H\varphi)v+\varphi(\Delta v)-2\langle d\varphi,dv\rangle
\end{array}$$
So

$$Lv=qv+\Delta v-\left\langle \frac{d(\varphi^2)}{\varphi^2},dv\right\rangle.$$
But for a positive function $\mu$, we have:

$$\Delta_{\mu}v=-\frac{1}{\mu}div(\mu grad v)=\Delta v-\frac{1}{\mu}\langle d\mu,dv\rangle, $$
hence the result.

\cqfd
\textit{End of the proof of Theorem (\ref{main result}):}\\
Applying Theorem (\ref{potentiel support compact}) to $L$, we deduce that $L$ has a finite number of non-positive eigenvalues. Therefore the same is true for $H$.
\cqfd

\subsection{Regularity questions}

In this section, we consider the case of a non-smooth potential $V$. We show that what we have proved remains true under a milder regularity assumption on $V$:

\begin{thm}\label{V non-smooth}

Let $M$ be a complete Riemannian manifold. Let $V\in L_{loc}^\infty$ be a potential such that $H=\Delta+V$ is bounded from below.\\
If there exists a positive function $\varphi\in W_{loc}^{1,2}$ such that $H\varphi=0$ weakly outside a compact set, then\\ 
$Card(Spec(H)\bigcap (-\infty,0])$, the number of non-positive bound states of $H=\Delta+V$, is finite.
 
\end{thm}
\textit{Proof of Theorem (\ref{V non-smooth}) :}\\
We will use the following result (cf \cite{Gilbarg-Trudinger}, Theorem 8.34): 

\begin{lem}

Let $V\in L_{loc}^\infty$ be a potential, and $H:=\Delta+V$. Let $u\in W_{loc}^{1,2}$ satisfying $Hu=0$ weakly inside a smooth, open, subset $\Omega$.\\ 
Then for all $\alpha\in (0,1)$ and all $\Omega'\subset\subset\Omega$, $u\in C^{1,\alpha}(\Omega')$.

\end{lem}
Given this Lemma, we first explain that we can assume, by modifying $\varphi$ on a compact set, that $\varphi\in C^{1,\alpha}_{loc}$. Let $K$ be a compact subset such that $H\varphi=0$ outside $K$, take $\tilde{\Omega}$ an open set such that $\tilde{\Omega}\subset\subset M\setminus K$. Let $\rho\in C_0^\infty(M)$ be a cut-off function such that $\rho\equiv1$ on $\tilde{K}=M\setminus{\tilde{\Omega}}$. Define $u=\rho .1+(1-\rho)\varphi$; $u\in C^{1,\alpha}_{loc}\cap W^{1,2}$ and $u>0$. In addition, we have:

\begin{lem}\label{Hu=0}
 
As a distribution, $Hu\in L^\infty$, and $Hu=0$ outside a compact set.

\end{lem}
\textit{Proof of Lemma (\ref{Hu=0}):}\\
We have:

$$Hu=H(\rho)+H((1-\rho)\varphi),$$
so that, given the the fact that $\varphi\in C^{1,\alpha}_{loc}(\tilde{\Omega})$, it is enough to prove that the following formula holds in the sense of distributions:

$$H((1-\rho)\varphi)=(\Delta(1-\rho))\varphi+2\langle d\rho,d\varphi\rangle$$
Let $\psi\in C_0^\infty(M)$, then by definition

\begin{displaymath}
\begin{array}{rcl}
 \langle H((1-\rho)\varphi),\psi\rangle&=&\langle (1-\rho)\varphi,H\psi\rangle\\\\
 &=& \langle (1-\rho)\varphi,\Delta \psi\rangle+\langle (1-\rho)\varphi,V\psi\rangle
\end{array}\end{displaymath}
Since by hypothesis $H\varphi=0$ outside $K$, we can substract $0=\langle \varphi,H((1-\rho)\psi)\rangle$ to the right term. Furthermore,

$$\langle \varphi,H((1-\rho)\psi)\rangle=\langle \varphi,(\Delta(1-\rho))\psi\rangle+\langle \varphi,(1-\rho)\Delta\psi\rangle-2\langle \varphi d(1-\rho),d\psi\rangle+\langle(1-\rho)\varphi,V\psi\rangle.$$
Thus we get:

$$\begin{array}{ccc}
 \langle H((1-\rho)\varphi),\psi\rangle&=&-\langle\varphi,(\Delta(1-\rho))\psi\rangle+2\langle\varphi d(1-\rho),d\psi\rangle
\end{array}$$
Given that $\varphi$ is in $C^{1,\alpha}_{loc}$ outside $K$, we can integrate by parts: 

$$\langle\varphi d(1-\rho),d\psi\rangle=\langle d^*(\varphi d(1-\rho)),\psi\rangle,$$
and furthermore the usual formula:

$$d^*(\varphi d(1-\rho))=\varphi\Delta(1-\rho)-\langle d\varphi,d(1-\rho)\rangle$$
is valid. Hence

$$\langle H((1-\rho)\varphi),\psi\rangle=\langle\varphi,(\Delta(1-\rho))\psi\rangle + 2\langle \langle d\varphi,d\rho\rangle ,\psi\rangle,$$
which is the result.

\cqfd
Given Lemma (\ref{Hu=0}), we can assume that there is a positive function $\varphi\in C^{1,\alpha}_{loc}\cap W_{loc}^{1,2}$, satisfying $H\varphi=0$ outside a compact set and such that $H\varphi\in L^\infty$. We want to mimic the proof of Theorem (\ref{main result}), and for this purpose we must show that the result of Lemma (\ref{Doob}) still holds. The point here is that the computations in the proof of Lemma (\ref{Doob}) require to assume $\varphi\in C^2_{loc}$, but here we only have $\varphi\in C^{1,\alpha}_{loc}$. It is the aim of the next Lemma to overcome this difficulty:

\begin{lem}\label{Doob2}
 
For every $v\in C_0^\infty$, and $\varphi\in C^1_{loc}$,

$$H(\varphi v)=(H\varphi)v+\varphi(\Delta v)-2\langle d\varphi,dv\rangle$$
in the sense of distributions.

\end{lem}
\textit{Proof of Lemma (\ref{Doob2}):}
Let $\psi\in C_0^\infty(M)$. By definition ,

$$\begin{array}{rcl}
\langle H(\varphi v),\psi\rangle &=&\langle \varphi v,H\psi\rangle \\\\
&=&\langle \varphi,H(v\psi)\rangle-\langle \varphi,\psi\Delta v\rangle+2\langle \varphi dv,d \psi\rangle \\\\
&=& \langle vH(\varphi),\psi\rangle-\langle \varphi,\psi\Delta v\rangle +2\langle d^*(\varphi dv),\psi\rangle
\end{array}$$
Since $\varphi\in C^{1,\alpha}$, we have the formula

$$d^*(\varphi dv)=\varphi \Delta v-\langle d\varphi,dv\rangle,$$
so we get

$$\langle H(\varphi v),\psi\rangle=\langle vH(\varphi),\psi\rangle+\langle \varphi \Delta v,\psi\rangle-2\langle \langle d\varphi,dv\rangle,\psi\rangle,$$
whence the result.

\cqfd

Therefore, letting $L:=\varphi^{-1}H\varphi$, for every $v\in C_0^\infty$ we have 

$$Lv=\left(\Delta_{\varphi^2}+q\right)v,$$
where $q=\varphi^{-1}H\varphi$. We thus get the equality $L=\Delta_{\varphi^2}+q$ as operators on distributions.\\
If the potential $q$ is in $L^\infty$, then the proof of Theorem (\ref{main result}) works. But this is just a consequence of the fact that we have assumed by Lemma (\ref{Hu=0}) that $H\varphi\in L^\infty$ and $\varphi$ continuous.

\cqfd

\subsection{Alternative proofs of the main result:}

\paragraph{First alternative method} In this paragraph, we explain how to get the result of Theorem (\ref{potentiel support compact}), without the statement on the kernel, by a different method. Let $H:=L+V$. We first introduce some notations.\\\\
Let us denote by $N_\lambda(H)$ the cardinal of $Spec(H)\cap(-\infty,\lambda)$, i.e. $\sup\{dim(W)\}$, $W$ subspace $C_0^\infty$ on which the quadratic form $q-\lambda$ is negative (recall that $q$ is the quadratic form associated to $H$). 
For $K\subset M$ the closure of a smooth, relatively compact set in $M$, we denote by $N_{K,\lambda}$ (resp. $N_{M\setminus K,\lambda}$) the cardinal of $Spec(H_K)\cap(-\infty,\lambda)$ (resp. $Spec(H_{M\setminus K})\cap(-\infty,\lambda)$), where $H_K$ (resp. $H_{M\setminus K}$) is $H$ on $K$ (resp. $M\setminus K$) with Neumann boundary conditions. Equivalently, $N_{K,\lambda}=\sup\{dim(W)\}$ (resp. $N_{M\setminus K,\lambda}=\sup\{dim(W)\}$), where $W$ is a subspace of $C^{\infty}(K)$ (resp. of $C_0^\infty(M\setminus K)$) on which the quadratic form $q-\lambda$ (resp. $H_{M\setminus K}-\lambda$) is negative. With these notations, we have the following relatively classical result (see \cite{Reed-Simon4}, Chapter 15, although it is not stated as such):

\begin{lem}\label{Nb vp}
For all $\lambda\in \mathbb{R}$,

$$N_\lambda\leq N_{K,\lambda}+N_{M\setminus K,\lambda}$$
\end{lem}
\textit{Proof of Lemma (\ref{Nb vp}):}\\
Let $W$ a subspace of $C_0^\infty$ on which the quadratic form associated to $H-\lambda$ is negative. Let $\varphi\in W\setminus\{0\}$. We have:

$$q(\varphi)<\lambda ||\varphi||_2^2.$$
Let $\varphi_1=\varphi|_{K}$ and $\varphi_2=\varphi|_{M\setminus K}$. Then $\varphi_1\in C^\infty(K)$ and $\varphi_2\in C_0^\infty(M\setminus K)$. If we can show that either $\frac{q(\varphi_1)}{||\varphi_1||_2^2}<\lambda$ or $\frac{q(\varphi_2)}{||\varphi_2||_2^2}<\lambda$, then we have the result. Suppose it is not the case, then:

$$q(\varphi_1)\geq \lambda ||\varphi_1||_2^2,$$
and
$$q(\varphi_2)\geq \lambda ||\varphi_2||_2^2.$$

But $\varphi=\varphi_1+\varphi_2$, and since $\varphi_1$ and $\varphi_2$ have the intersection of their support of measure zero, $q(\varphi)=q(\varphi_1)+q(\varphi_2)$ and $||\varphi||_2^2=||\varphi_1||_2^2+||\varphi_2||_2^2$. Therefore, we obtain

$$q(\varphi)\leq \lambda||\varphi||_2^2,$$
which is a contradiction.

\cqfd
Now, by standard elliptic theory, $N_{K,\lambda}$ is finite, for all $K$ as above and all $\lambda\in\mathbb{R}$. Thus, in order to prove Theorem (\ref{potentiel support compact}), we only need to find some $K$ such that $N_{M\setminus K,0}$ is finite. Take $K$ smooth containing the support of $V$. Then $H_{M\setminus K}$ is simply $L$ with Neumann boundary conditions on $\partial K$. But $L$ is a non-negative operator, therefore $N_{M\setminus K,0}=0$.

\cqfd

\begin{rem}
 
At first sight, we could think that Lemma (\ref{Nb vp}), combined with Lemma (\ref{Fischer-Schoen}) would give a proof of Theorem (\ref{main result}) in the general case (i.e. if we do not assume $V$ to be compactly supported), but the issue here is that Lemma (\ref{Fischer-Schoen}) only gives the non-negativity of $H$ restricted to $M\setminus K$ with Dirichlet boundary conditions, and not Neumann boundary conditions. In general, the infimum of the spectrum for Dirichlet boundary conditions is greater than the infimum of the spectrum for Neumann ones, so the Neumann operator is not necessarily non-negative. Of course, if the solution $\varphi$ satisfies Neumann boundary conditions, then a small variation of Lemma (\ref{Fischer-Schoen}) shows that the Neumann operator is non-negative, but the method of Fischer-Colbrie \cite{Fischer} does not easily yield the existence of such a $\varphi$ under the assumption that the Morse index is finite.

\end{rem}

\paragraph{Second method} We explain how we can adapt what we have done so that it is not necessary to use the Doob transform argument at the end. We start by assumption with a positive function $\varphi$ which satisfies $H\varphi=0$ outside a compact set. We can assume that $\varphi$ is defined on $M$. Let $\tilde{V}$ be a non-negative, compactly supported potential in $L^\infty$, such that $\tilde{V}\geq |H\varphi|$. Then

$$(H+\tilde{V})\varphi\geq 0,$$
which implies by Theorem (\ref{Fischer-Schoen}) that $H+\tilde{V}\geq 0$. We let $L:=H+\tilde{V}$, which is a non-negative operator of Schrödinger type:

$$L=\Delta_\mu+W,\,W\in L^\infty_{loc},$$
but the main difference is that $W$ is not assumed to be non-negative anymore. We owe to Yehuda Pinchover the remark that the theory of parabolicity has also been developped in this context. Let us quote the main result of \cite{Pinchover-Tintarev}:

\begin{thm}

Denote by $q$ the quadratic form associated to $L$. We have the following dichotomy:

\begin{enumerate}
	\item Either $L$ has a weighted spectral gap, i.e. there exists a positive function $\chi$ such that for all $u\in C_0^\infty(M)$,
	
	$$\int_Mu^2\chi\leq q(u).$$
	In this case, we say that $L$ is non-parabolic, or subcritical. Furthermore, $L$ has Green functions.

	\item Or there is a non-negative, compactly supported function $\psi$ and a positive function $\chi$ such that for all $u\in C_0^\infty(M)$,
	
	$$\int_Mu^2\chi\leq q(u)+\left(\int_M\psi u\right)^2.$$
	In this case, we say that $L$ is parabolic, or critical.
\end{enumerate}

\end{thm}
By Cauchy-Schwarz, 

$$\left(\int_M\psi u\right)^2\leq \left(\int_M \psi\right)\left(\int_M\psi u^2\right)=C\int_M\psi u^2,$$
so that, we can assume, by increasing $\tilde{V}$ if necessary, that $L$ is non-parabolic. We are therefore in the situation where $L$ has a weighted spectral gap: there exists a positive function $\chi$ such that for all $u\in C_0^\infty(M)$,
	
$$\int_Mu^2\chi\leq q(u).$$
This implies that the caracterisation of non-parabolicity of Theorem (\ref{non-parabolicity}) holds true for $L$, even if its potential is not non-negative. This allows us to make the theory of section \ref{consequence parabolicite} work for $L$, and as a consequence to prove the following Lemma:

\begin{lem}

For every compactly supported potential $R$ in $L^p$, $\frac{n}{2}<p\leq\infty$, $L+R$ has a finite number (counting multiplicities) of non-positive eigenvalues.

\end{lem}
Since $L$ and $H$ differ by a compactly supported potential $\tilde{V}$, this yields the result for $H$.
\cqfd

\noindent\textit{We would like to thank G. Carron whose questions have motivated this work and for his dedication and patience, and P. Castillon and Y. Pinchover for suggesting us the alternative proofs of the main theorem.}

\bibliographystyle{plain}

\bibliography{bibliographie}

\end{document}